\documentclass[10pt]{article}
\usepackage{geometry}
\newgeometry{vmargin={15mm}, hmargin={17mm,17mm}}

\usepackage[utf8]{inputenc}
\usepackage{amsfonts}
\usepackage{amsthm}
\usepackage{algorithm,algpseudocode}
\usepackage{amsmath}
\usepackage{babel}
\usepackage{verbatim}
\usepackage{tikz}
\usepackage{subcaption}
\usepackage{enumitem}

\newcommand{\NN}{\mathbb{N}}

\newcommand{\cT}{{\cal T}}
\newcommand{\cG}{{\cal G}}
\newcommand{\cO}{{\cal O}}

\newtheorem{definition}{Definition}
\newtheorem{lemma}{Lemma}
\newtheorem{theorem}{Theorem}
\newtheorem{corollary}{Corollary}
\newtheorem{observation}{Observation}
\newtheorem{example}{Example}
\newtheorem*{openProblem}{Open problem}

\newcommand{\footremember}[2]{%
	\footnote{#2}
	\newcounter{#1}
	\setcounter{#1}{\value{footnote}}%
}

\title{Counting traversing Hamiltonian cycles in tiled graphs}

\author{
	Alen Vegi Kalamar\footremember{affil1}{Department of Mathematics and Computer Science, University of Maribor, Maribor, Slovenia}\footremember{affil2}{Comtrade Gaming, Maribor, Slovenia}
}

\date{}

\newcommand\hlight[1]{\tikz[overlay, remember picture,baseline=-\the\dimexpr\fontdimen22\textfont2\relax]\node[rectangle,fill=blue!50,rounded corners,fill opacity = 0.2,draw,thick,text opacity =1] {$#1$};}

\providecommand{\keywords}[1]{\textbf{Keywords:} #1}

\begin{document}
	
\maketitle

\begin{abstract}
    In this paper we extend counting of traversing Hamiltonian cycles from 2-tiled graphs to generalized tiled graphs. We further show that, for a fixed finite set of tiles, counting traversing Hamiltonian cycles can be done in linear time with respect to the size of such graph, implying counting Hamiltonian cycles in tiled graphs is fixed-parameter tractable.\\
    \\
\end{abstract}

\section{Introduction}

Interest in the enumeration and generation of Hamiltonian cycles in various families of graphs can be found in different domains of science.
In chemistry, theoretical physics, and byophysics research is connected to modelling of polymers \cite{bib:Cloizeaux1990}, polymer melting, protein folding \cite{bib:Kloczkowski1998} and the study of magnetic systems with $\mathcal{O}(n)$ symmetry \cite{bib:Jacobsen2007}.
In theory of algorithms, it is one of the most studied counting problems coming from the mathematical nanosciences \cite{bib:Montoya2021}.
In engineering and bioinformatics, the research is connected to security and intellectual property protection \cite{bib:Nishat2019} and path planning problems for robots and machine tools \cite{bib:Liang2020}.

Studies on the topic of counting Hamiltonian cycles in different families of graphs and the use of similar (matrix) approaches are not negligible. In 1990, a~characterization of Hamiltonian cycles of the Cartesian product $P_4 \square P_n$ was established~\cite{bib:Tosic1990}. In 1994, Kwong and Rogers developed a matrix method for counting Hamiltonian cycles in $P_m\square P_n$, obtaining exact results for \mbox{$m=4,5$ \cite{bib:Kwong1994}}. Their method was extended to arbitrarily large grids by Bodroža-\mbox{Pantić~et~al.\ \cite{bib:Bodroza1994}} and by Stoyan and \mbox{Strehl~\cite{bib:Stoyan}}. Later, Bodroža-Pantić~et~al.\ gave some explicit generating functions for the number of Hamiltonian cycles in graphs $P_m \square P_n$, $C_m \square P_n$ and $P_m \square C_n$ \cite{bib:Bodroza2013, bib:Bodroza2019, bib:Bodroza2022}. Recently, Đokić~et~al.\ presented two algorithms for determination of the number of (all spanning unions of cycles) 2-factors in three classes of grid graphs: the thin cylinder $TnC_m(n)$, torus $TG^{(p)}_m(n)$ and Klein bottle $KB^{(p)}_m(n)$, all of which had a width of $m$ \cite{bib:Dokic2023}. On the other hand, Vegi Kalamar et al.\ characterized all Hamiltonian cycles in 2-tiled graphs and efficiently counted them for a fixed family of tiles \cite{bib:VegiKalamar2021}.

In the present contribution, we build on results of \cite{bib:VegiKalamar2021}. The introduced $k$-traversing Hamiltonian cycles in tiled graphs are a generalization of zigzagging (1-traversing) and traversing (2-traversing) Hamiltonian cycles from 2-tiled graphs. We extend results to counting traversing Hamiltonian cycles in arbitrary family of tiled graphs.

We organize the remainder of the paper as follows. In Section \ref{sec:tiledGraphsTraversingHamiltonianCycles}, we define tiled graphs as a generalization of 2-tiled graphs, define traversing Hamiltonian cycles in tiled graphs and prove their possible existence. In Section \ref{sec:countingTraversingHamiltonianCycles}, we present an algorithm to count this type of Hamiltonian cycles and represent its value in a closed formula. In Section \ref{sec:openProblems}, we propose open problems for further research.

\section{Tiled graphs and traversing Hamiltonian cycles}

\label{sec:tiledGraphsTraversingHamiltonianCycles}

In this section, we introduce the concept of a tile as presented in \cite{bib:VegiKalamar2021} and extend it to a definition of a tiled graph. We define the concept of $k$-traversing Hamiltonian cycles in tiled graphs and prove their possible existence for certain values of $k$.

\begin{definition}
    \label{def:tile}
    A \textit{tile} is a triple $T=(G,x,y)$, consisting of a connected graph G and two sequences $x = (x_1, x_2, \ldots, x_{k})$ (\textit{left wall}) and $y = (y_1, y_2, \ldots, y_{l})$ (\textit{right wall}) of distinct vertices of $G$, with~no vertex of $G$ appearing in both $x$ and $y$. We call $T$ a $(k,l)$\textit{-tile}.
\end{definition}

\begin{definition}
    \label{def:tiledGraph}
    \hspace{0cm}
    \begin{enumerate}
        \item The tiles $T=(G,x,y)$ and $T'=(G',x',y')$ are \textit{compatible} whenever $|y| = |x'|.$
        \item A sequence $\cT = (T_0, T_1,\ldots, T_m)$ of tiles is \textit{compatible} if, for~each $i \in \{1,2,\ldots,m\}$, $T_{i-1}$ is compatible with $T_{i}$.
        \item The \textit{join} of compatible tiles $(G,x,y)$ and $(G',x',y')$ is the tile $T = (G,x,y) \otimes (G',x',y')$ for which the graph is obtained from disjoint union of $G$ and $G'$ by identifying the sequence $y$ term by term with the sequence $x'$.
        \item The \textit{join} of a compatible sequence $\cT = (T_0, T_1,\ldots, T_m)$ of tiles is defined as $\otimes \cT = T_0 \otimes T_1 \otimes \cdots \otimes T_m$.
        \item A tile $T$ is \textit{cyclically compatible} if $T$ is compatible with itself.
        \item For a cyclically compatible tile $T=(G,x,y)$, the~\textit{cyclization} of $T$ is the graph $\circ T$ obtained by identifying the respective vertices of $x$ with $y$.
        \item A cyclization of a cyclically compatible sequence of tiles $\cT$ is defined as $\circ \cT = \circ (\otimes \cT)$.
        \item A \textit{tiled} graph is a cyclization of a sequence of at least 3 tiles.
    \end{enumerate}
\end{definition}

The following lemma is crucial in understanding of a structure of Hamiltonian cycles in tiled graphs. It is an extension of claims 1-3 of Lemma~1 from \cite{bib:VegiKalamar2021} to tiled graphs.

\begin{lemma}
    \label{lm:cycleProperty}
    Let $C$ be a Hamiltonian cycle in a tiled graph $G = \circ (T_0, T_1, \ldots, T_m)$. 
    Then, we have the following:
    \begin{enumerate}
        \item \label{cycleProperty.1} $C = \bigcup\limits_{i=0}^{m} (C \cap T_{i})$.
        \item \label{cycleProperty.2} $C \cap T_{i}$ is a union of paths and isolated vertices.
        \item \label{cycleProperty.3} Let $v$ be a vertex of a component of $C \cap T_i$. Then, $v$ has degree 2 in $C \cap T_i$ or~$v$ is a wall vertex.
    \end{enumerate}
\end{lemma}

\begin{proof}
    \hspace{0cm}
    \label{pr:cycleProperty}
    \begin{enumerate}
        \item $C = C \cap G = C \cap \big( \bigcup\limits_{i=0}^{m} T_{i} \big) \stackrel{\text{\tiny{Distributive law}}}{=} \bigcup\limits_{i=0}^{m} (C \cap T_{i}).$
        
        \item Let $K$ be a component of $C \cap T_i$. As~$C$ is a cycle, $K$ is a connected subgraph of $C$. Then, $K$ is either equal to $C$, a~path, or~a vertex. If~$K = C$, then $T_i$ contains all the vertices of $G$, a~contradiction to $m \geq 2$ (in at least one tile, $C$ does not contain all the vertices). The~claim~follows.
        
        \item Let $v$ be a vertex of $C \cap T_i$ of degree different from 2. As the maximum degree of vertex in $C$ is 2, $v$ has degree 1 or 0. If~$v$ is an internal vertex of $T_i$, its degree in $C = \bigcup\limits_{i=0}^{m} (C \cap T_{i})$ is equal to its degree in $C \cap T_i$. This contradicts $C$ being a cycle and the claim~follows.
    \end{enumerate}
\end{proof}

\begin{definition}
    \label{def:traversingCycle}
    Let $G = \circ (T_0, T_1, \ldots, T_m)$ be a tiled graph and $N_{i}$ the set of all isolated vertices in $C \cap T_i$. For $k \in \NN$, Hamiltonian cycle $C$ is $k$\textit{-traversing}, if $\forall i \in \{0, 1, \ldots, m\}$, $(C \cap T_{i}) \setminus N_i$ is a set of $k$ paths that start in a vertex of left wall and end in a vertex of right wall that cover all internal vertices of $T_i$.
\end{definition}

\begin{definition}
    Let $G = \circ (T_0, T_1, \ldots, T_m)$ be a tiled graph, where $\forall i, T_i$ is a $(k_i, l_i)$-tile. With $\min_{ws}(G)$ we denote the number $\min\{k_0, k_1, \ldots, k_m\}$. 
\end{definition}

The following result (using claim 4 of Lemma~1 from \cite{bib:VegiKalamar2021} as a basis) directly follows from claims of Lemma~\ref{lm:cycleProperty}.

\begin{corollary}
    \label{cr:traversingCycleExistence}
    Let $G = \circ (T_0, T_1, \ldots, T_m)$ be a tiled graph, where $\forall i, T_i$ is a $(k_i, l_i)$-tile. If a $k$-traversing Hamiltonian cycle exists in $G$, then $k \leq \min_{ws}(G)$.
\end{corollary}

\begin{proof}
    By Claim~\ref{cycleProperty.3} in Lemma~\ref{lm:cycleProperty}, paths start and end in a wall vertex. Each traversing path is a distinct non-degenerate path and has at least two unique wall vertices, one from left wall and one from right wall. For a tile $T_i$, the number of such paths is at most $\min\{k_i, l_i\}$ and the claim for the whole graph follows from the fact that $G = \circ (T_0, T_1, \ldots, T_m)$.
\end{proof}

\section{Counting traversing Hamiltonian cycles in tiled graphs}

\label{sec:countingTraversingHamiltonianCycles}

In this section, we introduce structures that enable us to describe an algorithm for counting traversing Hamiltonian cycles and prove its correctness. We further show that in the case of a fixed set of tiles, this algorithm is efficient.

Let $G = \circ (T_0, T_1, \ldots, T_m)$ be a tiled graph, where $\forall i$, $T_i$ is a $(k_i, l_i)$-tile.
For a tile $T_i$, let $k \in [\min\{k_i, l_i\}]$.
Let ${}_kL^e_i$ be a set of strings of left endvertices of length $k$ (out of $k_i$ options) and ${}_kR^e_i$ be a set of strings of right endvertices of length $k$ (out of $l_i$ options), both ordered lexicographically from smallest to largest. Then $|_kL^e_i| = \binom{k_i}{k}k!$ and $|_kR^e_i| = \binom{l_i}{k}k!$. Those sets will be used to model left-to-right transition modes.
Let $\epsilon$ denote an empty string. Additionally, let ${}_kL^r_i$ be a set of all binary strings of length $k_i - k$ (for remaining $k_i - k$ left wall vertices) and ${}_kR^r_i$ a set of all binary strings of length $l_i - k$ (for remaining $l_i - k$ left wall vertices),
$$
    _kL^r_i = 
    \begin{cases}
        \{ \epsilon \}, &\quad k = k_i\\
        \{0,1\}^{k_i - k}, &\quad\text{otherwise}
    \end{cases},
    \quad\quad
    _kR^r_i = 
    \begin{cases}
        \{ \epsilon \}, &\quad k = l_i\\
        \{0,1\}^{l_i - k}, &\quad\text{otherwise}
    \end{cases},
$$
both ordered lexicographically from largest to smallest. Then $|_kR^e_i| = 2^{k_i - k}$ and $|_kL^r_i| = 2^{l_i - k}$. Those sets will be used to model the coverage of wall vertices with paths within this (bit value 1) or a neighboring tile (bit value 0). For element $x$ from sets ${}_kL^r_i$ and ${}_kR^r_i$, let $\overline{x}$ denote a string where all bits are flipped and $\overline{\epsilon} = \epsilon$.

Let ${}^{i}_kA$ be a $|_kL^e_i| \times |_kR^e_i|$ block matrix
$$
    ^{i}_kA =
    \begin{bmatrix}
        ^{i}_kA_{b_1, c_1} & ^{i}_kA_{b_1, c_2} & \ldots & ^{i}_kA_{b_1, c_{|_kR^e_i|}} \\
        ^{i}_kA_{b_2, c_1} & ^{i}_kA_{b_2, c_2} & \ldots & ^{i}_kA_{b_2, c_{|_kR^e_i|}} \\
        \vdots & \vdots & \ddots & \vdots \\
        ^{i}_kA_{b_{|_kL^e_i|}, c_1} & ^{i}_kA_{b_{|_kL^e_i|}, c_2} & \ldots & ^{i}_kA_{b_{|_kL^e_i|}, c_{|_kR^e_i|}} \\
    \end{bmatrix},
$$
where, $\forall p \in [|_kL^e_i|]$, $b_p \in {}_kL^e_i$ ($p$-th element in ${}_kL^e_i$) and $\forall r \in [|_kR^e_i|]$, $c_r \in {}_kR^e_i$ ($r$-th element in ${}_kR^e_i$), ${}^{i}_kA_{b_p, c_r}$ is a $|_kL^r_i| \times |_kR^r_i|$ matrix of form
$$
    \begin{bmatrix}
        ^{i}_ka_{b_p, c_r}^{d_1, e_1} & ^{i}_ka_{b_p, c_r}^{d_1, e_2} & \ldots & ^{i}_ka_{b_p, c_r}^{d_1, e_{|_kR^r_i|}} \\
        ^{i}_ka_{b_p, c_r}^{d_2, e_1} & ^{i}_ka_{b_p, c_r}^{d_2, e_2} & \ldots & ^{i}_ka_{b_p, c_r}^{d_2, e_{|_kR^r_i|}} \\
        \vdots & \vdots & \ddots & \vdots \\
        ^{i}_ka_{b_p, c_r}^{d_{|_kL^r_i|}, e_1} & ^{i}_ka_{b_p, c_r}^{d_{|_kL^r_i|}, e_2} & \ldots & ^{i}_ka_{b_p, c_r}^{d_{|_kL^r_i|}, e_{|_kR^r_i|}} \\
    \end{bmatrix},
$$
where, $\forall u \in [|_kL^r_i|]$, $d_u \in {}_kL^r_i$ ($u$-th element in ${}_kL^r_i$), and $\forall v \in [|_kR^r_i|]$, $e_v \in {}_kR^r_i$ ($v$-th element in ${}_kR^r_i$), ${}^{i}_ka_{b_p, c_r}^{d_u, e_v}$ represents the number of possible combinations in $T_i$ for $k$ distinct paths covering all internal vertices that have ordered left wall endvertices $b_p$ and ordered right wall endvertices $c_r$ (defining paths $b_p^1 \to c_r^1, b_p^2 \to c_r^2, \ldots, b_p^k \to c_r^k$), $d_u$ describing if remaining $k_i - k$ left wall vertices are/are not part of any $k$ paths and $e_v$ describing if remaining $l_i - k$ right wall vertices are/are not part of any $k$ paths.

\tikzstyle{node}=[circle, draw, fill=black,inner sep=0pt, minimum width=4pt]

\begin{example}
    \label{ex:example1}
    Let $T$ be a (2,3)-tiled graph. Then there may exist 1-traversing and 2-traversing paths:
    \begin{itemize}[leftmargin=5pt]
        \item[] \textbf{1-traversing:}
        In this case, the sets $L^{e}, R^{e}, L^{r}, R^{r}$ are the following:
        \begin{align*}
            _1L^{e} &= \{1, 2\},\\
            _1R^{e} &= \{1, 2, 3\},\\
            _1L^{r} &= \{1, 0\},\\
            _1R^{r} &= \{11, 10, 01, 00\}.
        \end{align*}
        Then
        $$
            _1A =
            \begin{bmatrix}
                _1A_{1,1} & _1A_{1,2} & _1A_{1,3} \\
                _1A_{2,1} & _1A_{2,2} & _1A_{2,3}
            \end{bmatrix},
        $$
        where $\forall i \in L^{e}, \forall j \in R^{e}$
        $$
            _1A_{i,j} =
            \begin{bmatrix}
                _1a_{i,j}^{1,11} & _1a_{i,j}^{1,10} & _1a_{i,j}^{1,01} & _1a_{i,j}^{1,00} \\
                _1a_{i,j}^{0,11} & _1a_{i,j}^{0,10} & _1a_{i,j}^{0,01} & _1a_{i,j}^{0,00} \\
            \end{bmatrix}.
        $$
        For example, ${}_1a_{1,2}^{1, 10}$ denotes the number of possible combinations in $T$ for one path covering all internal vertices, starting in $x_1$, ending in $y_2$, where the first remaining left wall vertex ($x_2$) is part of the path and the first remaining right wall vertex ($y_1$) is part of the path, but the second one ($y_3$) is not.

        \begin{figure}[H]
            \centering
            \begin{subfigure}{0.49\textwidth}
                \centering
                \begin{tikzpicture}[scale=0.7]
        			\node[node] at (0, 1) {};
                    \node[label={$x_2$}] at (-0.5,0.6) {};
                    
        			\node[node] at (0, 3) {};
                    \node[label={$x_1$}] at (-0.5,2.6) {};
                    \node[circle,draw=red] at (0,3) {};
                    
        			\node[node] at (4, 0) {};
                    \node[label={$y_3$}] at (4.5,-0.4) {};
                    \draw[red] (3.8,-0.2) -- (4.2,0.2) (4.2,-0.2) -- (3.8,0.2);
                    
        			\node[node] at (4, 2) {};
                    \node[label={$y_2$}] at (4.5,1.6) {};
                    \node[circle,draw=red] at (4,2) {};
                    
        			\node[node] at (4, 4) {};
                    \node[label={$y_1$}] at (4.5,3.6) {};
                    
                    \draw[dotted, black] (0,3) -- (0,1) -- (4,4) -- (4,2);
    		      \end{tikzpicture}
            \end{subfigure}
            \hfill
            \begin{subfigure}{0.49\textwidth}
                \centering
        		\begin{tikzpicture}[scale=0.7]
        			\node[node] at (0, 1) {};
                    \node[label={$x_2$}] at (-0.5,0.6) {};
                    
        			\node[node] at (0, 3) {};
                    \node[label={$x_1$}] at (-0.5,2.6) {};
                    \node[circle,draw=red] at (0,3) {};
                    
        			\node[node] at (4, 0) {};
                    \node[label={$y_3$}] at (4.5,-0.4) {};
                    \draw[red] (3.8,-0.2) -- (4.2,0.2) (4.2,-0.2) -- (3.8,0.2);
                    
        			\node[node] at (4, 2) {};
                    \node[label={$y_2$}] at (4.5,1.6) {};
                    \node[circle,draw=red] at (4,2) {};
                    
        			\node[node] at (4, 4) {};
                    \node[label={$y_1$}] at (4.5,3.6) {};
                    
                    \draw[dotted, black] (0,3) -- (4,4) -- (0,1) -- (4,2);
    		      \end{tikzpicture}
            \end{subfigure}
    		\caption{Possible combinations for ${}_1a_{1,2}^{1, 10}$.}
    	\end{figure}
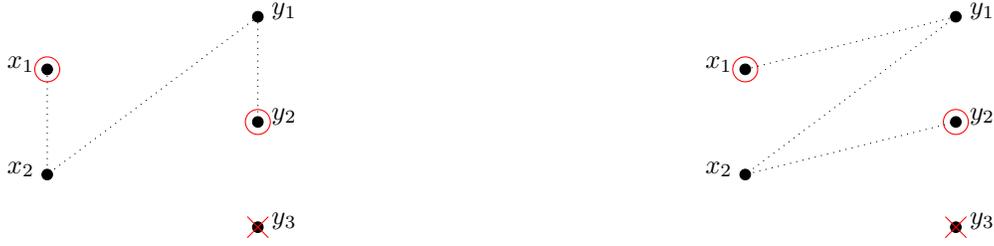
        
        \item[] \textbf{2-traversing:}
        In this case, the sets $L^{e}, R^{e}, L^{r}, R^{r}$ are the following:
        \begin{align*}
            _2L^{e} &= \{12, 21\},\\
            _2R^{e} &= \{12, 13, 21, 23, 31, 32\},\\
            _2L^{r} &= \{\epsilon\},\\
            _2R^{r} &= \{1, 0\}.
        \end{align*}
        Then
        $$
            _2A =
            \begin{bmatrix}
                _2A_{12,12} & _2A_{12,13} & _2A_{12,21} & _2A_{12,23} & _2A_{12,31} & _2A_{12,32} \\
                _2A_{21,12} & _2A_{21,13} & _2A_{21,21} & _2A_{21,23} & _2A_{21,31} & _2A_{21,32}
            \end{bmatrix},
        $$
        where $\forall i \in L^{e}, \forall j \in R^{e}$,
        $$
            _2A_{i,j} =
            \begin{bmatrix}
                _2a_{i,j}^{\epsilon,1} & _2a_{i,j}^{\epsilon,0}\\
            \end{bmatrix}.
        $$
        For example, ${}_2a_{21,23}^{\epsilon, 1}$ denotes the number of possible combinations in $T$ for two distinct paths covering all internal vertices, the first starting in $x_2$ and ending in $y_2$, and second one starting in $x_1$ and ending in $y_3$, where there is no remaining left wall vertex and the first remaining right wall vertex ($y_1$) is part of one of the paths.

        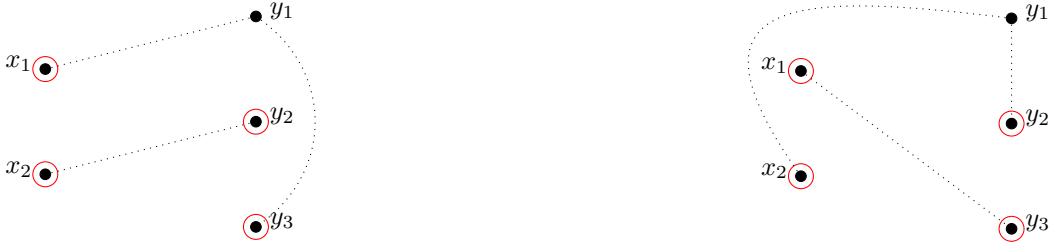
\begin{figure}[H]
            \centering
            \begin{subfigure}{0.49\textwidth}
                \centering
        		\begin{tikzpicture}[scale=0.7]
                    \node[node] at (0, 1) {};
                    \node[label={$x_2$}] at (-0.5,0.6) {};
                    \node[circle,draw=red] at (0,1) {};
                    
                    \node[node] at (0, 3) {};
                    \node[label={$x_1$}] at (-0.5,2.6) {};
                    \node[circle,draw=red] at (0,3) {};
                    
                    \node[node] at (4, 0) {};
                    \node[label={$y_3$}] at (4.5,-0.4) {};
                    \node[circle,draw=red] at (4,0) {};
                    
                    \node[node] at (4, 2) {};
                    \node[label={$y_2$}] at (4.5,1.6) {};
                    \node[circle,draw=red] at (4,2) {};
                    
                    \node[node] at (4, 4) {};
                    \node[label={$y_1$}] at (4.5,3.6) {};
                    
                    \draw[dotted, black] (0,1) -- (4,2) (0,3) -- (4,4);
                    \draw[dotted, black] (4,4) .. controls (5.5,3) and (5.5,1) .. (4,0);
                  \end{tikzpicture}
            \end{subfigure}
            \hfill
            \begin{subfigure}{0.49\textwidth}
                \centering
                \begin{tikzpicture}[scale=0.7]
                    \node[node] at (0, 1) {};
                    \node[label={$x_2$}] at (-0.5,0.6) {};
                    \node[circle,draw=red] at (0,1) {};
                    
                    \node[node] at (0, 3) {};
                    \node[label={$x_1$}] at (-0.5,2.6) {};
                    \node[circle,draw=red] at (0,3) {};
                    
                    \node[node] at (4, 0) {};
                    \node[label={$y_3$}] at (4.5,-0.4) {};
                    \node[circle,draw=red] at (4,0) {};
                    
                    \node[node] at (4, 2) {};
                    \node[label={$y_2$}] at (4.5,1.6) {};
                    \node[circle,draw=red] at (4,2) {};
                    
                    \node[node] at (4, 4) {};
                    \node[label={$y_1$}] at (4.5,3.6) {};
                    
                    \draw[dotted, black] (0,3) -- (4,0) (4,4) -- (4,2);
                    \draw[dotted, black] (0,1) .. controls (-2.5,4.5) and (0,4.5) .. (4,4);
                  \end{tikzpicture}
            \end{subfigure}
    		\caption{Possible combinations for ${}_2a_{21,23}^{\epsilon, 1}$.}
    	\end{figure}
    \end{itemize}
\end{example}

Let ${}^{i}_k\overline{A}$ be a $|_kL^e_i| \times |_kR^e_i|$ block matrix
$$
    ^{i}_k\overline{A} =
    \begin{bmatrix}
        ^{i}_k\overline{A}_{b_1, c_1} & ^{i}_k\overline{A}_{b_1, c_2} & \ldots & ^{i}_k\overline{A}_{b_1, c_{|_kR^e_i|}} \\
        ^{i}_k\overline{A}_{b_2, c_1} & ^{i}_k\overline{A}_{b_2, c_2} & \ldots & ^{i}_k\overline{A}_{b_2, c_{|_kR^e_i|}} \\
        \vdots & \vdots & \ddots & \vdots \\
        ^{i}_k\overline{A}_{b_{|_kL^e_i|}, c_1} & ^{i}_k\overline{A}_{b_{|_kL^e_i|}, c_2} & \ldots & ^{i}_k\overline{A}_{b_{|_kL^e_i|}, c_{|_kR^e_i|}} \\
    \end{bmatrix},
$$

where, $\forall p \in [|_kL^e_i|]$, $b_p \in {}_kL^e_i$ ($p$-th element in ${}_kL^e_i$) and $\forall r \in [|_kR^e_i|]$, $c_r \in {}_kR^e_i$ ($r$-th element in ${}_kR^e_i$), ${}^{i}_k\overline{A}_{b_p, c_r}$ is a $|_kL^r_i| \times |_kR^r_i|$ matrix of form
$$
    \begin{bmatrix}
        ^{i}_ka_{b_p, c_r}^{\overline{d_1}, e_1} & ^{i}_ka_{b_p, c_r}^{\overline{d_1}, e_2} & \ldots & ^{i}_ka_{b_p, c_r}^{\overline{d_1}, e_{|_kR^r_i|}} \\
        ^{i}_ka_{b_p, c_r}^{\overline{d_2}, e_1} & ^{i}_ka_{b_p, c_r}^{\overline{d_2}, e_2} & \ldots & ^{i}_ka_{b_p, c_r}^{\overline{d_2}, e_{|_kR^r_i|}} \\
        \vdots & \vdots & \ddots & \vdots \\
        ^{i}_ka_{b_p, c_r}^{\overline{d_{|_kL^r_i|}}, e_1} & ^{i}_ka_{b_p, c_r}^{\overline{d_{|_kL^r_i|}}, e_2} & \ldots & ^{i}_ka_{b_p, c_r}^{\overline{d_{|_kL^r_i|}}, e_{|_kR^r_i|}} \\
    \end{bmatrix}.
$$

\begin{example}
    Using tile from Example~\ref{ex:example1}, we check the ${}_k\overline{A}$ matrices:
    \begin{itemize}[leftmargin=5pt]
        \item[] \textbf{1-traversing:}
        $$
            _1\overline{A} =
            \begin{bmatrix}
                _1\overline{A}_{1,1} & _1\overline{A}_{1,2} & _1\overline{A}_{1,3}\\
                _1\overline{A}_{2,1} & _1\overline{A}_{2,2} & _1\overline{A}_{2,3}
            \end{bmatrix},
        $$
        where $\forall i \in L^{e}, \forall j \in R^{e}$
        $$
            _1\overline{A}_{i,j} =
            \begin{bmatrix}
                _1a_{i,j}^{\overline{1},11} & _1a_{i,j}^{\overline{1},10} & _1a_{i,j}^{\overline{1},01} & _1a_{i,j}^{\overline{1},00} \\
                _1a_{i,j}^{\overline{0},11} & _1a_{i,j}^{\overline{0},10} & _1a_{i,j}^{\overline{0},01} & _1a_{i,j}^{\overline{0},00} \\
            \end{bmatrix}
            =
            \begin{bmatrix}
                _1a_{i,j}^{0,11} & _1a_{i,j}^{0,10} & _1a_{i,j}^{0,01} & _1a_{i,j}^{0,00} \\
                _1a_{i,j}^{1,11} & _1a_{i,j}^{1,10} & _1a_{i,j}^{1,01} & _1a_{i,j}^{1,00} \\
            \end{bmatrix}.
        $$
        
        \item[] \textbf{2-traversing:}
        $$
            _2\overline{A} =
            \begin{bmatrix}
                _2\overline{A}_{12,12} & _2\overline{A}_{12,13} & _2\overline{A}_{12,21} & _2\overline{A}_{12,23} & _2\overline{A}_{12,31} & _2\overline{A}_{12,32} \\
                _2\overline{A}_{21,12} & _2\overline{A}_{21,13} & _2\overline{A}_{21,21} & _2\overline{A}_{21,23} & _2\overline{A}_{21,31} & _2\overline{A}_{21,32}
            \end{bmatrix},
        $$
        where $\forall i \in L^{e}, \forall j \in R^{e}$
        $$
            _2\overline{A}_{i,j} =
            \begin{bmatrix}
                _2a_{i,j}^{\overline{\epsilon},1} & _2a_{i,j}^{\overline{\epsilon},0}\\
            \end{bmatrix}
            =
            \begin{bmatrix}
                _2a_{i,j}^{\epsilon,1} & _2a_{i,j}^{\epsilon,0}\\
            \end{bmatrix}.
        $$
    \end{itemize}
\end{example}

\begin{definition}
    Let $T_i$ be a $(k_i, l_i)$-tile and $T_{i+1}$ a $(k_{i+1}, l_{i+1})$-tile that are compatible ($l_i = k_{i+1}$). For $k \in [\min\{k_i, l_i, l_{i+1}\}]$, let ${}_k^{i}A$ be a matrix that belongs to $T_i$ and ${}_k^{i+1}\overline{A}$ a matrix that belongs to $T_{i+1}$. For a tile $T_i \otimes T_{i+1}$, we define a matrix ${}_k^{i, i+1}A$ as
    $$
        ^{i, i + 1}_k A = {}^{i}_k A \cdot {}^{i+1}_k \overline{A},
    $$
    where $\cdot$ represents the block matrix multiplication
    $$
        ^{i, i + 1}_kA_{b_p, c_r} = \sum\limits_{x \in {}_kR^e_i = {}_kL^e_{i+1}} \Big({}^{i}_k A_{b_p, x} \cdot {}^{i + 1}_k \overline{A}_{x, c_r}\Big).
    $$
\end{definition}

\begin{lemma}
    Matrix ${}_k^{i, i+1}A$ is a $|_kL^e_i| \times |_kR^e_{i+1}|$ block matrix of form
    $$
        \begin{bmatrix}
        ^{i, i+1}_kA_{b_1, c_1} & ^{i, i+1}_kA_{b_1, c_2} & \ldots & ^{i, i+1}_kA_{b_1, c_{|_kR^e_i|}} \\
        ^{i, i+1}_kA_{b_2, c_1} & ^{i, i+1}_kA_{b_2, c_2} & \ldots & ^{i, i+1}_kA_{b_2, c_{|_kR^e_i|}} \\
        \vdots & \vdots & \ddots & \vdots \\
        ^{i, i+1}_kA_{b_{|_kL^e_i|}, c_1} & ^{i, i+1}_kA_{b_{|_kL^e_i|}, c_2} & \ldots & ^{i, i+1}_kA_{b_{|_kL^e_i|}, c_{|_kR^e_i|}} \\
    \end{bmatrix},
    $$
    where, $\forall p \in [|_kL^e_i|]$, $b_p \in {}_kL^e_i$ ($p$-th element in ${}_kL^e_i$) and $\forall r \in [|_kR^e_{i+1}|]$, $c_r \in {}_kR^e_{i+1}$ ($r$-th element in ${}_kR^e_{i+1}$), ${}^{i, i+1}_kA_{b_p, c_r}$ is a $|_kL^r_i| \times |_kR^r_{i+1}|$ matrix of form
    $$
        \begin{bmatrix}
            ^{i, i+1}_ka_{b_p, c_r}^{d_1, e_1} & ^{i, i+1}_ka_{b_p, c_r}^{d_1, e_2} & \ldots & ^{i, i+1}_ka_{b_p, c_r}^{d_1, e_{|_kR^r_{i+1}|}} \\
            ^{i, i+1}_ka_{b_p, c_r}^{d_2, e_1} & ^{i, i+1}_ka_{b_p, c_r}^{d_2, e_2} & \ldots & ^{i, i+1}_ka_{b_p, c_r}^{d_2, e_{|_kR^r_{i+1}|}} \\
            \vdots & \vdots & \ddots & \vdots \\
            ^{i, i+1}_ka_{b_p, c_r}^{d_{|_kL^r_i|}, e_1} & ^{i, i+1}_ka_{b_p, c_r}^{d_{|_kL^r_i|}, e_2} & \ldots & ^{i, i+1}_ka_{b_p, c_r}^{d_{|_kL^r_i|}, e_{|_kR^r_{i+1}|}} \\
        \end{bmatrix},
    $$
    where, $\forall u \in [|_kL^r_i|]$, $d_u \in {}_kL^r_i$ ($u$-th element in ${}_kL^r_i$), and $\forall v \in [|_kR^r_{i+1}|]$, $e_v \in {}_kR^r_{i+1}$ ($v$-th element in ${}_kR^r_{i+1}$), ${}^{i, i+1}_ka_{b_p, c_r}^{d_u, e_v}$ represents the number of possible combinations in $T_i \otimes T_{i+1}$ for $k$ distinct paths covering all internal vertices that have ordered left wall endvertices $b_p$ and ordered right wall endvertices $c_r$ (defining paths $b_p^1 \to c_r^1, b_p^2 \to c_r^2, \ldots, b_p^k \to c_r^k$), $d_u$ describing if remaining $k_i - k$ left wall vertices are/are not part of any $k$ paths and $e_v$ describing if remaining $l_{i+1} - k$ right wall vertices are/are not part of any $k$ paths.
\end{lemma}

\begin{proof}
    Because tiles $T_i$ and $T_{i+1}$ are compatible, ${}_kR^e_i = {}_kL^e_{i+1}$ and ${}_kR^r_i = {}_kL^r_{i+1}$.
    Since ${}_k^{i}A$ is a block matrix of dimension $|_kL^e_i| \times |_kR^e_i|$ and  ${}_k^{i+1}\overline{A}$ is a block matrix of dimension $|_kL^e_{i+1}| \times |_kR^e_{i+1}|$, by definition of block matrix multiplication, ${}_k^{i, i+1}A$ is a block matrix of dimension $|_kL^e_i| \times |_kR^e_{i+1}|$.
    Since each block in ${}_k^{i}A$ is of dimension $|_kL^r_i| \times |_kR^r_i|$ and each block in ${}_k^{i+1}\overline{A}$ of dimension $|_kL^r_{i+1}| \times |_kR^r_{i+1}|$, by definition of block matrix multiplication, each block ${}_k^{i, i+1}A_{b_p, c_r}$ in ${}_k^{i,i+1}A$ is of dimension $|_kL^r_i| \times |_kR^r_{i+1}|$.
    
    By definition of block matrix multiplication, element ${}^{i, i+1}_ka_{b_p, c_r}^{d_u, e_v}$ of block ${}_k^{i, i+1}A_{b_p, c_r}$ is calculated as
    $$
        ^{i, i+1}_ka_{b_p, c_r}^{d_u, e_v} = \sum\limits_{x \in {}_kR^e_i}\Big( \sum\limits_{y \in {}_kR^r_i} ({}^{i}_ka_{b_p, x}^{d_u, y} \cdot {}^{i+1}_ka_{x, c_r}^{\overline{y}, e_v} ) \Big).
    $$
    We have to show that $\forall x \in {}_kR^e_i, \forall y \in {}_kR^r_i$,
    $$
        ^{i}_ka_{b_p, x}^{d_u, y} \cdot {}^{i+1}_ka_{x, c_r}^{\overline{y}, e_v}
    $$
    generates described paths in $T_i \otimes T_{i+1}$.
    
    Let $x \in {}_kR^e_i$ and $y \in {}_kR^r_i$. Then ${}^{i}_ka_{b_p, x}^{d_u, y}$ represents the number of possible combinations in $T_i$ for $k$ distinct paths covering all internal vertices that have ordered left wall endvertices $b_p$ and ordered right wall endvertices $x$ (defining paths $b_p^1 \to x^1, b_p^2 \to x^2, \ldots, b_p^k \to x^k$), $d_u$ describing if remaining $k_i - k$ left wall vertices are/are not part of any $k$ paths and $y$ describing if remaining $l_i - k$ right wall vertices are/are not part of any $k$ paths.
    On the other side, ${}^{i+1}_ka_{x, c_r}^{\overline{y}, e_v}$ represents the number of possible combinations in $T_{i+1}$ for $k$ distinct paths covering all internal vertices that have ordered left wall endvertices $x$ and ordered right wall endvertices $c_r$ (defining paths $x^1 \to c_r^1, x^2 \to c_r^2, \ldots, x^k \to c_r^k$), $\overline{y}$ describing if remaining $k_{i+1} - k$ left wall vertices are/are not part of any $k$ paths and $e_v$ describing if remaining $l_{i+1} - k$ right wall vertices are/are not part of any $k$ paths.
    
    It is clear that if we identify those $k$ paths by endvertices $x$ (defining paths $b_p^1 \to x^1 \to c_r^1, b_p^2 \to x^2 \to c_r^2, \ldots, b_p^k \to x^k \to c_r^k$), we cover all internal vertices in $T_i$ and $T_{i+1}$, but also remaining vertices on shared wall (paths in $T_i$ cover remaining right wall vertices in $T_i$, for which bit value in $y$ is 1, and paths in $T_{i+1}$ cover remaining left wall vertices in $T_{i+1}$, for which bit value in $\overline{y}$ is 1) and the result follows.
\end{proof}

\begin{example}
    We check out the product ${}_k^1A \cdot {}_k^2\overline{A}$, where $T_1$ is a tile from Example~\ref{ex:example1} and $T_2$ the tile that we get if we switch walls in $T_1$ (those tiles are compatible).
    \begin{itemize}[leftmargin=5pt]
        \item[] \textbf{1-traversing:}
        \begin{align*}
            ^{1,2}_1A
            &=
            {}^1_1A \cdot {}^2_1\overline{A}\\
            &=
            \begin{bmatrix}
                ^1_1A_{1,1} & ^1_1A_{1,2} & ^1_1A_{1,3} \\
                ^1_1A_{2,1} & ^1_1A_{2,2} & ^1_1A_{2,3}
            \end{bmatrix}
            \cdot
            \begin{bmatrix}
                ^2_1\overline{A}_{1,1} & ^2_1\overline{A}_{1,2}\\
                ^2_1\overline{A}_{2,1} & ^2_1\overline{A}_{2,2}\\
                ^2_1\overline{A}_{3,1} & ^2_1\overline{A}_{3,2}
            \end{bmatrix}\\
            &=
            \begin{bmatrix}
                _1^{1,2}A_{1,1} & _1^{1,2}A_{1,2} \\
                _1^{1,2}A_{2,1} & _1^{1,2}A_{2,2}
            \end{bmatrix},
        \end{align*}
        where $\forall i \in L^{e}, \forall j \in R^{e}$
        $$
            ^{1,2}_1A_{i,j} = \sum\limits_{x \in {}_kR^e} {}^1_1A_{i,x} \cdot {}^2_1\overline{A}_{x,j}.
        $$
        
        \item[] \textbf{2-traversing:}
        \begin{align*}
            ^{1,2}_2A
            &=
            {}^1_2A \cdot {}^2_2\overline{A}\\
            &=
            \begin{bmatrix}
                ^1_2A_{12,12} & ^1_2A_{12,13} & ^1_2A_{12,21} & ^1_2A_{12,23} & ^1_2A_{12,31} & ^1_2A_{12,32} \\
                ^1_2A_{21,12} & ^1_2A_{21,13} & ^1_2A_{21,21} & ^1_2A_{21,23} & ^1_2A_{21,31} & ^1_2A_{21,32}
            \end{bmatrix}
            \cdot
            \begin{bmatrix}
                ^2_2\overline{A}_{12,12} & ^2_2\overline{A}_{12,21}\\
                ^2_2\overline{A}_{13,12} & ^2_2\overline{A}_{13,21}\\
                ^2_2\overline{A}_{21,12} & ^2_2\overline{A}_{21,21}\\
                ^2_2\overline{A}_{23,12} & ^2_2\overline{A}_{23,21}\\
                ^2_2\overline{A}_{31,12} & ^2_2\overline{A}_{31,21}\\
                ^2_2\overline{A}_{32,12} & ^2_2\overline{A}_{32,21}\\
            \end{bmatrix}\\
            &=
            \begin{bmatrix}
                _2^{1,2}A_{12,12} & _2^{1,2}A_{12,21} \\
                _2^{1,2}A_{21,12} & _2^{1,2}A_{21,21}
            \end{bmatrix},
        \end{align*}
        where $\forall i \in L^{e}, \forall j \in R^{e}$
        $$
            ^{1,2}_2A_{i,j} = \sum\limits_{x \in {}_kR^e} {}^1_2A_{i,x} \cdot {}^2_2\overline{A}_{x,j}.
        $$
    \end{itemize}
\end{example}

\begin{definition}
    For $j > 1$, we define
    $$
        ^{i, i + j}_k A = {}^{i, i+j-1}_k A \cdot {}^{i+j}_k \overline{A}.
    $$
\end{definition}

\begin{definition}
    For $b = b^1 b^2 \ldots b^k \in L^e_0$ and $z \in [k]$, let
    $$
        b_{\hookrightarrow z} = b^{k-z+1}\ldots b^{k} b^{1} \ldots b^{k-z}
    $$
    be the $z$-shift of $b$.
\end{definition}

\begin{theorem}
    Let $G = \circ (T_0, T_1, \ldots, T_m)$ be a tiled graph and let $k \leq \min_{ws}(G)$. The number of distinct $k$-traversing Hamiltonian cycles in $G$ is equal to
    $$
    {}_kTHC(G) =
    \begin{cases}
        \sum\limits_{b \in {}_kL^e_0}\sum\limits_{s = 1}^{|_kL^r_0|} {}^{0,m}_ka_{b, b}^{d_s, d_{|_kL^r_0| + 1 - s}}, &\quad \text{if } k = 1,\\
        \frac{1}{k!}\sum\limits_{b \in {}_kL^e_0}\sum\limits_{z = 1}^{k-1}\sum\limits_{s = 1}^{|_kL^r_0|} {}^{0,m}_ka_{b, b_{\hookrightarrow z}}^{d_s, d_{|_kL^r_0| + 1 - s}}, &\quad \text{if } k > 1.
    \end{cases}
    $$
\end{theorem}

\begin{proof} 
    Let ${}^{0,m}_kA$ be a matrix that belongs to $T_0 \otimes T_1 \otimes \cdots \otimes T_m$. Then it is a $|_kL^e_0| \times |_kR^e_m|$ block matrix of form
    $$
        \begin{bmatrix}
            ^{0,m}_kA_{b_1, c_1} & ^{0,m}_kA_{b_1, c_2} & \ldots & ^{0,m}_kA_{b_1, c_{|_kR^e_m|}} \\
            ^{0,m}_kA_{b_2, c_1} & ^{0,m}_kA_{b_2, c_2} & \ldots & ^{0,m}_kA_{b_2, c_{|_kR^e_m|}} \\
            \vdots & \vdots & \ddots & \vdots \\
            ^{0,m}_kA_{b_{|_kL^e_0|}, c_1} & ^{0,m}_kA_{b_{|_kL^e_0|}, c_2} & \ldots & ^{0,m}_kA_{b_{|_kL^e_0|}, c_{|_kR^e_m|}}
        \end{bmatrix},
    $$
    where $\forall p \in [|_kL^e_0|]$, $b_p \in {}_kL^e_0$ ($p$-th element in ${}_kL^e_0$), $\forall r \in [|_kR^e_m|]$, $c_r \in$ ($r$-th element in ${}_kR^e_m$), ${}^{0,m}_kA_{b_p, c_r}$ is a $|_kL^r_0| \times |_kR^r_m|$ matrix of form
    $$
        \begin{bmatrix}
            ^{0,m}_ka_{b_p, c_r}^{d_1, e_1} & ^{0,m}_ka_{b_p, c_r}^{d_1, e_2} & \ldots & ^{0,m}_ka_{b_p, c_r}^{d_1, e_{|_kR^r_m|}} \\
            ^{0,m}_ka_{b_p, c_r}^{d_2, e_1} & ^{0,m}_ka_{b_p, c_r}^{d_2, e_2} & \ldots & ^{0,m}_ka_{b_p, c_r}^{d_2, e_{|_kR^r_.|}} \\
            \vdots & \vdots & \ddots & \vdots \\
            ^{0,m}_ka_{b_p, c_r}^{d_{|_kL^r_0|}, e_1} & ^{0,m}_ka_{b_p, c_r}^{d_{|_kL^r_0|}, e_2} & \ldots & ^{0,m}_ka_{b_p, c_r}^{d_{|_kL^r_0|}, e_{|_kR^r_m|}}
        \end{bmatrix}.
    $$
    
    Since graph $G$ is a result of identifications of left and right wall vertices in tile $T_0 \otimes T_1 \otimes \cdots \otimes T_m$, $k_0 = l_m$ and so ${}_kL^e_0 = {}_kR^e_m$ and ${}_kL^r_0 = {}_kR^r_m$. We get that matrix ${}^{0,m}_kA$ is a $|_kL^e_0| \times |_kL^e_0|$ block matrix of form
    $$
        \begin{bmatrix}
            ^{0,m}_kA_{b_1, b_1} & ^{0,m}_kA_{b_1, b_2} & \ldots & ^{0,m}_kA_{b_1, b_{|_kL^e_0|}} \\
            ^{0,m}_kA_{b_2, b_1} & ^{0,m}_kA_{b_2, b_2} & \ldots & ^{0,m}_kA_{b_2, b_{|_kL^e_0|}} \\
            \vdots & \vdots & \ddots & \vdots \\
            ^{0,m}_kA_{b_{|_kL^e_0|}, b_1} & ^{0,m}_kA_{b_{|_kL^e_0|}, b_2} & \ldots & ^{0,m}_kA_{b_{|_kL^e_0|}, b_{|_kL^e_0|}}
        \end{bmatrix},
    $$
    where $\forall p, r \in [|_kL^e_0|]$, $b_p, b_r \in {}_kL^e_0$ ($p$-th and $r$-th element in ${}_kL^e_0$), ${}^{0,m}_kA_{b_p, b_r}$ is a $|_kL^r_0| \times |_kL^r_0|$ matrix of form
    $$
        \begin{bmatrix}
            ^{0,m}_ka_{b_p, b_r}^{d_1, d_1} & ^{0,m}_ka_{b_p, b_r}^{d_1, d_2} & \ldots & ^{0,m}_ka_{b_p, b_r}^{d_1, d_{|_kL^r_0|}} \\
            ^{0,m}_ka_{b_p, b_r}^{d_2, d_1} & ^{0,m}_ka_{b_p, b_r}^{d_2, d_2} & \ldots & ^{0,m}_ka_{b_p, b_r}^{d_2, d_{|_kL^r_0|}} \\
            \vdots & \vdots & \ddots & \vdots \\
            ^{0,m}_ka_{b_p, b_r}^{d_{|_kL^r_0|}, d_1} & ^{0,m}_ka_{b_p, b_r}^{d_{|_kL^r_0|}, d_2} & \ldots & ^{0,m}_ka_{b_p, b_r}^{d_{|_kL^r_0|}, d_{|_kL^r_0|}}
        \end{bmatrix},
    $$
    where, $\forall u, v \in [|_kL^r_0|]$, $d_u, d_v \in {}_kL^r_0$ ($u$-th and $v$-th element in ${}_kL^r_0$), ${}^{i}_ka_{b_p, b_r}^{d_u, d_v}$ represents the number of possible combinations in $T_0 \otimes T_1 \otimes \cdots \otimes T_m$ for $k$ distinct paths covering all internal vertices that have ordered left wall endvertices $b_p$ and ordered right wall endvertices $b_r$ (defining paths $b_p^1 \to b_r^1, b_p^2 \to b_r^2, \ldots, b_p^k \to b_r^k$), $d_u$ describing remaining $k_0 - k$ left wall vertices are/are not part of any $k$ paths and $d_v$ describing remaining $k_0 - k$ right wall vertices are/are not part of any $k$ paths.
    
    To get the number ${}_kTHC(G)$ from the matrix ${}^{0,m}_kA$, we need to take into account only block matrices ${}^{0,m}_kA_{b_p, b_r}$ for which mapping of consecutive elements of $b_p$ into consecutive elements of $b_r$ ($b_p^1 \mapsto b_r^1, \ldots, b_p^k \mapsto b_r^k$) is a permutation of a set $\{b_p^1, \ldots, b_p^k\}$ that can be decomposed into one disjoint cycle. Otherwise each sequence of paths identified by a disjoint cycle in the decomposition can generate a separate cycle. For $k = 1$, $\forall b \in {}_kL^e_0$, $b \to b$ is the desired permutation, and, for $k > 1$, it is easy to see, that $\forall b \in {}_kL^e_0$, mappings $b \to b_{\hookrightarrow z}$, where $z \in [k-1]$, are the desired permutations.
    
    From such block matrices ${}^{0,m}_kA_{b, c}$, only entries ${}^{0,m}_ka_{b, c}^{d_s, \overline{d_s}}$ add to the number ${}_kTHC(G)$, since they cover all remaining vertices on the shared wall between tiles $T_m$ and $T_0$ exactly once.
    Because of the introduced lexicographic order on sets ${}_kL^r_i, {}_kR^r_i$,
    \begin{align*}
        \forall s \in [|_kL^r_0|]: \overline{d_s} &= d_{|_kL^r_0| + 1 - s}.
    \end{align*}
    Hence, if $k = 1$, $\forall b \in {}_kL^e_0$,
    $$
        {}^{0,m}_ka_{b, b}^{d_s, \overline{d_s}} = {}^{0,m}_ka_{b, b}^{d_s, d_{|_kL^r_0| + 1 - s}},
    $$
    and, if $k > 1$, $\forall b \in {}_kL^e_0, \forall z \in [k-1],$
    $$
        {}^{0,m}_ka_{b, b_{\hookrightarrow z}}^{d_s, \overline{d_s}} = {}^{0,m}_ka_{b, b_{\hookrightarrow z}}^{d_s, d_{|_kL^r_0| + 1 - s}}
    $$
    are values that contribute to ${}_kTHC(G)$.
    
    However, since there are $k!$ block matrices that calculate same cycles (ones that generate permutations over the same set), the total sum must be divided by $k!$.
\end{proof}

\begin{observation}
    \hspace{0cm}
    \begin{enumerate}
        \item If $k = 1$, for $b \in {}_kL^e_0$, set
        $$
            \Bigg\{ {}^{0,m}_ka_{b, b}^{d_s, \overline{d_s}}\ \Bigg| \ s \in [|_kL^r_0|] \Bigg\}
        $$
        forms an anti-diagonal of block matrix ${}^{0,m}_kA_{b, b}$.
        \item If $k > 1$, for $b \in {}_kL^e_0$ and $z \in [k-1]$, set
        $$
            \Bigg\{ {}^{0,m}_ka_{b, b_{\hookrightarrow z}}^{d_s, \overline{d_s}}\ \Bigg| \ s \in [|_kL^r_0|] \Bigg\} 
        $$
        forms an anti-diagonal of block matrix ${}^{0,m}_kA_{b, b_{\hookrightarrow z}}$.
    \end{enumerate}
\end{observation}

Using the above algorithm, some results from \cite{bib:VegiKalamar2021} can be recreated. In it, 2-tiled graphs were investigated.
\begin{example}[Counting traversing Hamiltonian cycles in 2-tiled graphs]
    \label{ex:traversingIn2TiledGraphs}
    {\ } \\ 
    For 2-tiled graphs, 1-traversing (denoted as zigzagging in \cite{bib:VegiKalamar2021}) and 2-traversing (denoted as traversing in \cite{bib:VegiKalamar2021}) Hamiltonian cycles may exist:
    \begin{itemize}[leftmargin=5pt]
        \item[] \textbf{1-traversing:}
        The matrix ${}^{0,m}_kA$ is a $2 \times 2$ block matrix, where each block is of dimension $2 \times 2$:
        $$
            \left[\begin{array}{c|c}
                ^{0,m}_kA_{1, 1} & ^{0,m}_kA_{1, 2}\\
                \hline
                ^{0,m}_kA_{2, 1} & ^{0,m}_kA_{2, 2}\\
            \end{array}\right]
            =
            \left[\begin{array}{cc|cc}
                ^{0,m}_ka_{1, 1}^{1, 1} & ^{0,m}_ka_{1, 1}^{1, 0} & ^{0,m}_ka_{1, 2}^{1, 1} & ^{0,m}_ka_{1, 2}^{1, 0}\\
                ^{0,m}_ka_{1, 1}^{0, 1} & ^{0,m}_ka_{1, 1}^{0, 0} & ^{0,m}_ka_{1, 2}^{0, 1} & ^{0,m}_ka_{1, 2}^{0, 0}\\
                \hline
                ^{0,m}_ka_{2, 1}^{1, 1} & ^{0,m}_ka_{2, 1}^{1, 0} & ^{0,m}_ka_{2, 2}^{1, 1} & ^{0,m}_ka_{2, 2}^{1, 0}\\
                ^{0,m}_ka_{2, 1}^{0, 1} & ^{0,m}_ka_{2, 1}^{0, 0} & ^{0,m}_ka_{2, 2}^{0, 1} & ^{0,m}_ka_{2, 2}^{0, 0}
            \end{array}\right].
        $$
        Then
        \begin{align*}
            {}_1THC(G) &= {}^{0,m}_ka_{1, 1}^{1, 0} + {}^{0,m}_ka_{1, 1}^{0, 1} + {}^{0,m}_ka_{2, 2}^{1, 0} + {}^{0,m}_ka_{2, 2}^{0, 1}.
        \end{align*}
        
        \item[] \textbf{2-traversing:}
        The matrix ${}^{0,m}_kA$ is a $2 \times 2$ block matrix, where each block is of dimension $1 \times 1$:
        $$
            \left[\begin{array}{c|c}
                ^{0,m}_kA_{12, 12} & ^{0,m}_kA_{12, 21}\\
                \hline
                ^{0,m}_kA_{21, 12} & ^{0,m}_kA_{21, 21}\\
            \end{array}\right]
            =
            \left[\begin{array}{c|c}
                ^{0,m}_ka_{12, 12}^{\epsilon, \epsilon} & ^{0,m}_ka_{12, 21}^{\epsilon, \epsilon}\\
                \hline
                ^{0,m}_ka_{21, 12}^{\epsilon, \epsilon} & ^{0,m}_ka_{21, 21}^{\epsilon, \epsilon}
            \end{array}\right].
        $$
        Then
        $$
            {}_2THC(G) = \frac{1}{2!} ({}^{0,m}_ka_{12, 21}^{\epsilon, \epsilon} + ^{0,m}_ka_{21, 12}^{\epsilon, \epsilon}) = {}^{0,m}_ka_{12, 21}^{\epsilon, \epsilon}.
        $$
    \end{itemize}
\end{example}

\begin{corollary}
    \label{cr:traversingInTiledGraphs}
    Let $\cT$ be a finite family of tiles, and let $\cG$ be a family of cyclizations of finite sequence of such tiles. There exists and algorithm that yields, for each graph $G \in \cG, \forall k \leq \min_{ws}(G)$, the number ${}_kTHC(G)$. For a fixed set $\cT$, the running time of the algorithm is linear in the number of tiles (and hence vertices) of $G$.
\end{corollary}

\begin{proof}
    For a fixed family of tiles, we precalculate matrices ${}^{i}_kA$ (and so ${}^{i}_k\overline{A}$). The time complexity to compute the matrix ${}^{0,m}_kA$ is $\cO(m)$ and additional $\cO(1)$ is needed to get the searched sum from it.
\end{proof}

\section{Open problems}

\label{sec:openProblems}

In this paper, we presented a solution to count traversing Hamiltonian cycles in tiled graphs. The crucial part of doing it was the identification of matrices ${}^{i}_kA$. In Corollary~\ref{cr:traversingInTiledGraphs}, we used the fact that in a finite set of tiles, the generation of these matrices can be thought of as a constant problem. In general, this is not true, so we propose the following open problem:

\begin{openProblem}
    Find an (efficient) algorithm to generate block matrices ${}^{i}_kA$.
\end{openProblem}
    
Traversing Hamiltonian cycles are such that their intersection with the tile results in paths that start in one wall and end in another wall of the tile. It is obvious, however, that there are also types of Hamiltonian cycles such that their intersection with the tile also allows paths that start and end in the same wall of the tile (flanking Hamiltonian cycles as an example in 2-tiled graphs \cite{bib:VegiKalamar2021}). Therefore we additionally propose the following open problem:

\begin{openProblem}
    Identify and count other types of Hamiltonian cycles in tiled graphs.
\end{openProblem}

\section{Acknowledgments}

I am grateful to Drago Bokal for his critical comments that encouraged me to improve the presentation of the paper. This research was funded by Slovenian Research Agency (research project J1-2452).


\begin{thebibliography}{999}
	\bibitem{bib:Cloizeaux1990} des Cloizeaux, J.; Jannik, G. \textit{Polymers in Solution: Their Modelling and Structure}; Clarendon Press: Oxford, MS, USA, 1990.
 
    \bibitem{bib:Kloczkowski1998} Kloczkowski, A.; Jernigan, R.L. Transfer matrix method for enumeration and generation of compact self-avoiding walks. I. Square lattices. {\em J. Chem. Phys.} {\bf 1998}, {\em 109}, 5134–-5146.
    
    \bibitem{bib:Jacobsen2007} Jacobsen, J.L. Exact enumeration of Hamiltonian circuits, walks and chains in two and three dimensions. {\em J. Phys. A Math. Theor.} {\bf 2007}, {\em 40}, 14667-–14678.
    
    \bibitem{bib:Montoya2021} Montoya, J.A. On the Counting Complexity of Mathematical Nanosciences. {\em MATCH Commun. Math. Comput. Chem.} {\bf 2021}, {\em 86}, 453-–488.

    \bibitem{bib:Nishat2019} Nishat, R.I.; Whitesides, S. Reconfiguring Hamiltonian Cycles in L-Shaped Grid Graphs. {\em Graph-Theor. Concepts Comput. Sci.} {\bf 2019}, {\em 21}, 325-–337.
    
    \bibitem{bib:Liang2020} Liang, T.C.; Chakrabarty, K.; Karri, R. Programmable daisychaining of microelectrodes to secure bioassay IP in MEDA biochips. {\em IEEE Trans. Very Large Scale Integr. (VLSI) Syst.} {\bf 2020}, {\em 25}, 1269-–1282.

	\bibitem{bib:Tosic1990} Tošić, R.; Bodroža-Pantić, O.; Kwong, Y.H.H.; Straight, H.J. On the number of Hamiltonian cycles of $P_4 \times P_n$. {\em Indian J. Pure Appl. Math.} {\bf 1990}, {\em 21}, 403--409.
	
	\bibitem{bib:Kwong1994} Kwong, Y.H.H.; Rogers, D.G. A matrix method for counting Hamiltonian cycles on grid graphs. {\em Eur. J. Comb.} {\bf 1994}, {\em 15}, 277--283.
	
	\bibitem{bib:Bodroza1994} Bodroža-Pantić, O.; Tošič, R. On the number of 2-factors in rectangular lattice graphs. {\em Publ. De L'Institut Math.} {\bf 1994}, {\em 56}, 23--33.
	
	\bibitem{bib:Stoyan} Stoyan, R.; Strehl, V. Enumeration of Hamiltonian circuits in rectangular grids. {\em J. Comb. Math. Comb. Comput.} {\bf 1996}, {\em 21}, 109--127.
	
	\bibitem{bib:Bodroza2013} Bodroža-Pantić, O.; Pantić, B.; Pantić, I.; Bodroža-Solarov, M. Enumeration of Hamiltonian cycles in some grid graphs. {\em MATCH Commun. Math. Comput. Chem.} {\bf 2013}, {\em 70}, 181--204.
	
	\bibitem{bib:Bodroza2019} Bodroža-Pantić, O.; Kwong, H.; Doroslovački, R.; Pantić, M. Enumeration of Hamiltonian cycles on a thick grid cylinder---Part I: Non-contractible Hamiltonian cycles. {\em Appl. Anal. Discret. Math.} {\bf 2019}, {\em 13}, 28--60.
	
	\bibitem{bib:Bodroza2022} Bodroža-Pantić, O.; Kwong, H.; Dokić, J.; Doroslovački, R.; Pantić, M. Enumeration of Hamiltonian cycles on a thick grid cylinder---Part II: Contractible Hamiltonian cycles. {\em Appl. Anal. Discret. Math.} {\bf 2022}, {\em 16}, 246--287.

    \bibitem{bib:Dokic2023} Đokić, J.; Doroslovački, K.; Bodroža-Pantić, O. A Spanning Union of Cycles in Thin Cylinder, Torus and Klein Bottle Grid Graphs. {\em Mathematics} {\bf 2023}, {\em 11, 846}.
 
    \bibitem{bib:VegiKalamar2021} Vegi Kalamar, A.; Žerak, T.; Bokal, D. Counting Hamiltonian Cycles in 2-Tiled Graphs. {\em  Mathematics} {\bf 2021}, {\em 9(6), 693}, 1--27.
\end{thebibliography}
\end{document}